\documentstyle[12pt]{article}
\font\goth=eufm10 scaled 1400

\def\g{\gamma}

\def\s{\sigma}

\def\S1{\hbox{\rm S$^1$}}
\def\Vect{\hbox{\rm Vect}}

\def\pds#1,#2{\langle #1\mid #2\rangle} 
\def\f#1,#2,#3{#1\colon#2\to#3} 

\def\hfl#1{{\buildrel{#1}\over{\hbox to
12mm{\rightarrowfill}}}}

\chardef\s=110
\chardef\g=103

\begin{document}

\title{Deforming the Lie algebra of vector fields on $S^1$ inside
the Poisson algebra on $\dot T^*S^1$}

\author{V. Ovsienko {\small and}
C. Roger}

\date{}

\maketitle

{\abstract{We study deformations of the standard embedding of the
 Lie algebra $\Vect(S^1)$ of smooth vector fields on the circle, 
into the Lie algebra 
of functions on the cotangent bundle $T^*S^1$
(with respect to the Poisson bracket).
We consider two analogous but different problems:
(a) formal deformations of the standard embedding
of $\Vect(S^1)$ into the Lie algebra 
of functions on 
$\dot T^*S^1:=T^*S^1\!\setminus\!S^1$ 
which are Laurent polynomials on fibers, and
(b) polynomial deformations of the $\Vect(S^1)$ subalgebra
inside the Lie algebra of formal Laurent series on $\dot T^*S^1$.
}}

\vskip 7cm


\thispagestyle{empty}

\vfill\eject

\section{Introduction}

{\bf 1.1 The standard embedding}. 

\noindent
The Lie algebra $\Vect(M)$ of vector fields on a manifold $M$ has a
natural embedding
 into the Poisson Lie algebra of functions on $T^*M$.
It is defined by the standard action of the Lie algebra
of vector fields on the cotangent bundle. 
Using the local Darboux coordinates 
$(x,\xi)=(x^1,\dots,x^n,\xi_1,\dots,\xi_n)$ on $T^*M$,
the explicit formula is:
\begin{equation}
\pi(X)
=X\xi 
\label{pi}
\end{equation}
where $X$ is a vector field: $X=\sum_{i=1}^n X^i(x)\partial /\partial x^i$ and
$X\xi=\sum_{i=1}^n X^i(x)\xi_i$.

The main purpose of this paper is to study 
deformations of the standard embedding (\ref{pi}). 

\vskip 0,3cm

\noindent
{\bf 1.2 Deformations inside $C^{\infty}(T^*M)$}.

\noindent
Consider the Poisson Lie algebra of smooth functions on $T^*M$
for an orientable manifold $M$.
In this case, the problem of deformation of the embedding (\ref{pi})
has an elementary solution.
The $\Vect(M)$ embedding (\ref{pi}) 
into
$C^{\infty}(T^*M)$
has the unique (well-known) nontrivial deformation.
Indeed, given an arbitrary volume form on $M$,
the expression:
$$
\pi_{\lambda}(X)=
X\xi+\lambda{\rm div}X,
$$
where $\lambda\in\bf R$,
defines an embedding of $\Vect(M)$
into $C^{\infty}(T^*M)$.

The linear map:
$X\mapsto{\rm div}X$ 
is the unique nontrivial
1-cocycle on  $\Vect(M)$ with values in 
$C^{\infty}(M)\subset C^{\infty}(T^*M)$
(cf. \cite{FUK}).

\vskip 0,3cm

\noindent
{\bf 1.3 Two Poisson Lie algebras of formal symbols}.

\noindent
Let us consider
the following two Lie algebras of Poisson
on the cotangent bundle
with zero section removed:
$
\dot T^*M=T^*M\!\setminus\!M.
$

(a) The Lie algebra $A(M)$ of functions on $\dot T^*M$
which are
{\it Laurent polynomials} on fibers;

(b) The Lie algebra ${\cal A}(M)$ of {\it formal Laurent series} on $\dot T^*M$.

Lie algebras $A(M)$ and ${\cal A}(M)$ can be interpreted as
{\it classical limits of the algebra of formal symbols} of pseudo-differential
operators on
$M$. We will show that in this case one can expect much more
interesting results than those in the case of $C^{\infty}(M)$. 

In both cases, the Poisson bracket is defined by the usual formula:
$$
\{F,G\}=
\frac{\partial F}{\partial \xi}\frac{\partial G}{\partial x}-
\frac{\partial F}{\partial x}\frac{\partial G}{\partial \xi}.
$$

\section{Statement of the problem}

In this paper we will consider only the one-dimensional case: $M=S^1$ 
(analogous results hold for $M=\bf R$).

\vskip 0,3cm
\noindent
{\bf 2.1 Algebras $A(S^1)$ and ${\cal A}(S^1)$ in the one-dimensional case}.

\noindent
As vector spaces, Lie algebras 
$A(S^1)$ and ${\cal A}(S^1)$ have the following form:
$$
A(S^1):=C^{\infty}(S^1)\otimes{\bf C}[\xi,\xi^{-1}]
\;\;\;
{\rm and}
\;\;\;
{\cal A}(S^1):=C^{\infty}(S^1)\otimes{\bf C}[[\xi,\xi^{-1}]],
$$
where ${\bf C}[[\xi,\xi^{-1}]]$ is the space of Laurent series
in one formal indeterminate.

Elements of both algebras: 
$A(S^1)$ and ${\cal A}(S^1)$ can be written
in the following form:
$$
F(x,\xi,\xi^{-1})=\sum_{k\in{\bf Z}}\xi^kf_k(x),
$$
where the coefficients $f_k(x)$
are periodic functions: $f_k(x+2\pi)=f_k(x)$.
In the case of algebra $A(S^1)$, one supposes that
the coefficients $f_k\equiv0$, if $|k|$ is sufficiently large;
for ${\cal A}(S^1)$
the condition is: $f_k\equiv0$, if $k$ is sufficiently large.

\vskip 0,3cm
\noindent
{\bf 2.2 Formal deformations of $\Vect(S^1)$ inside $A(S^1)$}.

\noindent
We will study {\it one-parameter formal} deformations of the standard
embedding 
of $\Vect(S^1)$ into the Lie algebra
$A(S^1)$.
That means we study linear maps
$$
\pi^t:\Vect(S^1)
\to
A(S^1)[[t]]
$$
to the Lie algebra of series in
a formal parameter $t$.
Such a map has the following form:
\begin{equation}
\pi^t=
\pi+t\pi_1+t^2\pi_2+\cdots
\label{pit}
\end{equation}
where $\pi_k:\Vect(S^1)\to A(S^1)$
are some linear maps,
such that the formal homomorphism
condition is satisfied: 
$$
\pi_t([X,Y])=\{\pi_t(X),\pi_t(Y)\}.
$$
The general Nijenhuis--Richardson theory of formal deformations
of homomorphisms of Lie algebras will be discussed in the next
section.

\vskip 0,3cm

\noindent
{\bf 2.3 Polynomial deformations of the $\Vect(S^1)$ inside ${\cal A}(S^1)$}.

\noindent
We classify all the {\it polynomial} deformations of the standard
embedding (\ref{pi}) of
$\Vect(S^1)$ into ${\cal A}(S^1)$.
In other words, we consider homomorphisms 
of the following form:
\begin{equation}
\pi(c)=\pi+\sum_{k\in{\bf Z}}\pi_k(c)\xi^k
\label{pic}
\end{equation}
where $c=c_1,\dots,c_n\in\bf R$ (or $\bf C$) are parameters of deformations,
each linear map $\pi_k(c):\Vect(S^1)\to C^{\infty}(S^1)$ being polynomial in $c$,
$\pi_k(0)=0$ and $\pi_k\equiv0$ if $k>0$ is sufficiently large.

\vskip 0,3cm

\noindent
{\bf 3.4 Motivations}. 

\noindent
(a)
Lie algebras
of functions on a symplectic manifold have nontrivial formal deformations
linked with so-called {\it deformation quantization}.
The problem considered in this paper,
is original and have never been discussed in the literature.
However, this problem is inspired by
deformation quantization.

The geometric version of the problem,
deformations (up to symplectomorphism)
of zero section of the cotangent bundle 
$M\subset T^*M$,
has no nontrivial solutions.
Existence of nontrivial deformations in the algebraic
formulation that we consider here seems to be a manifestation of ``quantum
anomalies''.

Note, that interesting examples of deformations of 
Lie algebra homomorphisms related to deformation quantization
can be found in \cite{TAB}.

\vskip 0,3cm

\noindent
(b)
Lie algebras of vector fields and Lie algebras
of functions on a symplectic manifold, have both nice cohomology theories,
our idea is to link them together.

Lie algebras of vector fields have various nontrivial extensions.
The well-known example is the {\it Virasoro algebra}
defined as a central extension of $\Vect(S^1)$.
A series of nontrivial extensions of $\Vect(S^1)$
by modules of tensor-densities on $S^1$ were
constructed in  \cite{OR1},\cite{OR2}.
These extensions can be obtained,
using a (nonstandard) embedding of 
$\Vect(S^1)$ into $C^{\infty}(\dot T^*S^1)$,
by restriction of the deformation of $C^{\infty}(\dot T^*S^1)$ (see \cite{OR2}). 

We will show that deformations of the standard embedding
relate the Virasoro algebra to extensions of Poisson algebra on ${\bf T}^2$
defined by A.A. Kirillov (see \cite{KIR},\cite{ROG}).

\vskip 0,3cm

\noindent
(c)
The following quantum aspect of the considered problem: deformations of
embeddings of $\Vect(S^1)$ into the algebra of pseudodifferential 
operators on $S^1$, will be
treated in a subsequent article.

\section{Nijenhuis-Richardson theory}

Deformations of homomorphisms of  Lie algebras
were first considered in \cite{NR} (see also \cite{RIC}).
The Nijenhuis--Richardson theory is analogous
to the Gerstenhaber theory of formal deformations
of associative algebras (and Lie algebras) (see \cite{GER}),
related cohomological calculations are parallel.
Let us outline the main results of this theory.

\vskip 0,3cm

\noindent
{\bf 3.1 Equivalent deformations}.

\vskip 0,3cm

\noindent
{\bf Definition}. Two homomorphisms $\pi$ and $\pi'$ of a Lie algebra 
{\goth g}
to a Lie algebra 
{\goth h}
are {\it equivalent} (cf. \cite{NR}) if there exists an 
{\it interior} automorphism
$I$ of {\goth h} such that $\pi'=I\pi$.

\vskip 0,3cm

Let us specify this definition for the two problem formulated
in Sections 2.2 and 2.3.

\vskip 0,3cm

\noindent
(a) Two formal deformations (\ref{pit}) $\pi_t$ and $\pi'_t$ are equivalent
if there exists a linear map
$I_t:A(S^1)[[t]]\to A(S^1)[[t]]$ 
of the form:
$$
\matrix{
I_t=
{\rm exp}(t{\rm ad}_{F_1}+t^2{\rm ad}_{F_2}+\cdots)\hfill\cr\noalign{\smallskip}
\;\;\;\;={\rm id}+
t{\rm ad}_{F_1}+
t^2({\rm ad}^2_{F_1}/2+{\rm ad}_{F_2})+
\cdots\hfill\cr
}
$$
where $F_i\in A(S^1)$,
such that $\pi'_t=I_t\pi_t$.
It is natural to consider such an automorphism of $A(S^1)[[t]]$ 
as interior.

\vskip 0,3cm

\noindent
(b) An automorphism
$I(c):{\cal A}(S^1)\to{\cal A}(S^1)$ 
depending on the parameters
$c=c_1,\dots,c_n$,
which is of the following form:
$$
I(c)=
{\rm exp}(\sum_{i=1}^{n}c_i{\rm ad}_{F_i}+
c_ic_j{\rm ad}_{F_{ij}}+\cdots)
$$
where $F_i,F_{ij},\dots\in{\cal A}(S^1)$
is called interior.
Two polynomial deformations $\pi(c)$ and $\pi'(c)$ 
of the standard embedding
$\Vect(S^1)\hookrightarrow{\cal A}(S^1)$
are equivalent
if there exists an interior automorphism
$I(c)$,
such that $\pi'(c)=I(c)\pi(c)$.

\vskip 0,3cm

\noindent
{\bf 3.2 Infinitesimal deformations}.

\noindent
Deformations (\ref{pit}) and (\ref{pic}),
modulo second order terms in $t$ and $c$ respectively,
are called {\it infinitesimal}.
Infinitesimal deformations of a Lie algebra homomorphism
from {\goth g}
into 
{\goth h}
are classified by the first cohomology group
$H^1(${\goth g};{\goth h}$)$, {\goth h} being a 
{\goth g}-module through $\pi$.

Namely, the first order terms $\pi_1$ in (\ref{pit})
and $\left.\frac{\partial\pi(c)}{\partial c_i}\right|_{c=0}$ in (\ref{pic})
are {\it 1-cocycles}.
Two infinitesimal deformations are equivalent if and only if
the corresponding cocycles are cohomologous.

Conversely,
given a Lie algebra homomorphism
$\pi:${\goth g}
$\to$
{\goth h},
an arbitrary 1-cocycle 
$\pi_1\in Z^1(${\goth g};{\goth h}$)$
defines an infinitesimal deformation of $\pi$.

\vskip 0,3cm

\noindent
{\bf 3.3 Obstructions}. 

\noindent
The {\it integrability conditions} are conditions for existence of 
(formal or polynomial) deformation corresponding
to a given infinitesimal deformation.

\vskip 0,3cm

\noindent
(a)
The obstructions for existence of a formal deformation (\ref{pit})
belong to the second cohomology group
$H^2(${\goth g};{\goth h}$)$.
This follows from so-called deformation relation (see \cite{NR}):
\begin{equation}
d\pi_t+(1/2)[\pi_t,\pi_t]=O
\label{deq}
\end{equation}
where $[\pi_t,\pi_t]$
is a bilinear map from
{\goth g}
to
{\goth h}:
$$
[\pi_t,\pi_t](x,y):=\{\pi_t(x),\pi_t(y)\}+\{\pi_t(y),\pi_t(x)\}.
$$
Note that the deformation relation (\ref{deq}) is nothing but
a rewritten formal homomorphism relation (Section 1.4).

The equation (\ref{deq}) is equivalent to a series of nonlinear equations
concerning the maps $\pi_k$:
$$
d\pi_k=\sum_{i+j=k}[\pi_i,\pi_j].
$$
The right hand side of each equation is a 2-cocycle
and the equations have solutions if and only if the
corresponding cohomology classes vanish.

\vskip 0,3cm

\noindent
(b)
Analogous {\it necessary} conditions for existence
of a polynomial deformation (\ref{pic}) can be easily calculated.

\vskip 0,3cm

\noindent
{\bf 3.4 Remarks: polynomial deformations}.

\noindent
Deformations of algebraic structures
(as associative and Lie algebras, their modules and
homomorphisms)
{\it polynomially depending on parameters} are not very well studied.
There is no special version of the general theory adopted to this case
and the number of known examples is small
(see \cite{AMM}).

Theory of polynomial deformation seems to be richer
than those of formal ones.
The equivalence problem for polynomial deformation has additional interesting
aspects related to parameter transformations
(cf. Sections 5.4 and 5.5, formul{\ae} (\ref{par})).

\section{Polynomial deformations of the embedding of $\Vect(S^1)$
into the Lie algebra of formal Laurent series
on $T^*S^1$}

Consider the Poisson Lie algebra 
${\cal A}(S^1)$. 
The formula (\ref{pi})
defines an embedding of 
$\Vect(S^1)$ into this Lie algebra.

The following theorem is the main result of this paper.
It gives a classification of 
polynomial deformations of the subalgebra
$
\Vect(S^1)\subset {\cal A}(S^1).
$

\proclaim Theorem 1. Every nontrivial polynomial deformation of the standard
embedding of $\Vect(S^1)$ into 
${\cal A}(S^1)$
is equivalent to one of a two-parameter family of deformations 
given by the formula:
\begin{equation}
\pi_{\lambda,\mu}\bigg(f(x)\frac{d}{dx}\bigg)=
f\bigg(x+\frac{\lambda-\mu}{\xi}\bigg)\xi +
\mu f^{\prime}\bigg(x+\frac{\lambda-\mu}{\xi}\bigg)
\label{uni}
\end{equation}
where $\lambda,\mu\in{\bf R}$ or $\bf C$ are parameters of the
deformation; the expression in the right
hand side has to be interpreted as a formal (Taylor) series in $\xi$.
\par

A complete proof of this theorem is given in Sections 4 and 5.

\vskip 0,3cm

The explicit formula for the deformation $\pi_{\lambda,\mu}$ is as follows:
$$
\matrix{
\pi_{\lambda,\mu}(f(x)\frac{d}{dx})=f(x)\xi+\lambda f'(x)+
\left(\frac{\lambda^2}{\!2}-\frac{\mu^2}{\!2}\right)f''(x)\xi^{-1}
+\cdots\hfill\cr\noalign{\smallskip}
\;\;\;\;\;\;\;\;\;\;\;\;\;\;\;\;\;\;\;\;\;\;\;\;\;\;
\displaystyle+\left(\frac{\mu(\lambda-\mu)^k}{k!}+
\frac{\;(\lambda-\mu)^{k+1}}{\!\!\!\!(k+1)!}\right)
f^{(k+1)}(x)\xi^{-k}+\cdots\hfill\cr
}
\eqno(5')
$$

\vskip 0,3cm

\noindent
{\bf Remark}.
The formula (\ref{uni}) is a result of complicated calculations
which will be omitted.
We do not see any {\it a-priori} reason for its existence.

\vskip 0,3cm

To prove Theorem 1, we apply the Nijenhuis--Richardson theory.

\vskip 0,3cm

The first step is to classify infinitesimal deformations. 
One has to calculate the first cohomology of
$\Vect(S^1)$ with coefficients in ${\cal A}(S^1)$.
Then,
one needs the integrability condition under which an
infinitesimal deformation corresponds to
a polynomial one.

\vskip 0,3cm

\noindent
{\bf 4.1 Algebras ${\cal A}(S^1)$ and $A(S^1)$ as a $\Vect(S^1)$-modules}.

\noindent
Lie algebra $\Vect(S^1)$ is a subalgebra of ${\cal A}(S^1)$.
Therefore, ${\cal A}(S^1)$ is a $\Vect(S^1)$-module.

\vskip 0,3cm

\noindent
{\bf Definition}. Consider a 1-parameter family of $\Vect (S^1)$-actions on 
$C^{\infty}(S^1)$ given by
$$
L^{(\lambda)}_{f(x)\frac{d}{dx}}(a(x))=
f(x)a^{\prime}(x)-\lambda f^{\prime}(x)a(x)
$$
where $\lambda\in \bf R$.

\vskip 0,3cm

\noindent
Denote ${\cal F}_{\lambda}$ 
the $\Vect (S^1)$-module structure on $C^{\infty}(S^1)$
defined by this action.

\vskip 0,3cm

\noindent
{\bf Remark}. Geometricaly, $L^{(\lambda)}_{fd\!/\! dx}$
is the operator of Lie derivative on {\it tensor-densities}
of degree $-\lambda$. That means:
$a=a(x)(dx)^{-\lambda}.$

\proclaim Lemma 4.1. (i) The Lie algebra $A(S^1)$ is decomposed to
a direct sum of $\Vect (S^1)$-modules:
$$
A(S^1)=\oplus_{m\in{\bf Z}}{\cal F}_m.
$$
\hfill\break
(ii) The Lie algebra ${\cal A}(S^1)$
has the following decomposition
as a $\Vect (S^1)$-module:
$$
\displaystyle{\cal A}(S^1)=
\oplus_{m\geq0}{\cal F}_m\oplus
\Pi_{m<0}{\cal F}_m.
$$
\par

\vskip 0,3cm

\noindent
{\bf Proof}. Consider the subspace
of $A(S^1)$ and ${\cal A}(S^1)$ consisting of functions
of degree $m$ in $\xi$:
$a(x)\xi^m$.
This subspace
is a $\Vect (S^1)$-module isomorphic to ${\cal F}_m$.
One has:
$$
\{f(x)\xi,a(x)\xi^m\}=(fa^{\prime}-mf^{\prime}a)\xi^m
=L_{f\frac{d}{dx}}^{(m)}(a)\xi^m.
$$
Therefore, algebra of Laurent polynomials
is a direct sum of $\Vect (S^1)$-modules
${\cal F}_m$.

By definition, an element of algebra ${\cal A}(S^1)$ is
a formal series in $\xi$ 
with a finite number of terms of positive degree.

Lemma 4.1 is proven.

\vskip 0,3cm

\noindent
{\bf 4.2 Cohomology groups $H^1(\Vect(S^1);A(S^1))$
and $H^1(\Vect(S^1);{\cal A}(S^1))$}.

\noindent
It follows that the $\Vect(S^1)$-cohomology 
with coefficients in $A(S^1)$
is splitted into a direct sum:
$$
H^1(\Vect(S^1);A(S^1))=
\oplus_{m\in{\bf Z}}H^1(\Vect(S^1);{\cal F}_m).
$$
These cohomology groups are well known (see \cite{FUK}).
They are nontrivial if and only if $m=0,-1,-2$
and the corresponding
group of cohomology are one-dimensional.
Therefore, the space of first cohomology $H^1(\Vect(S^1);A(S^1))$ is
three-dimensional. 

It is clear that the same result holds for ${\cal A}(S^1)$:
$$
H^1(\Vect(S^1);A(S^1))
\;=\;
H^1(\Vect(S^1);{\cal A}(S^1))\;=\;{\bf R}^3
$$

The nontrivial cocycles generating the cohomology groups
$H^1(\Vect(S^1);\allowbreak A(S^1))$ and $H^1(\Vect(S^1);{\cal A}(S^1))$
are as follows:
$$
\matrix{
C_0(fd\!/\! dx)= f'\hfill\cr
C_1(fd\!/\! dx)= f''(dx)\hfill\cr
C_2(fd\!/\! dx)= f'''(dx)^2\hfill\cr
}
$$
with values in ${\cal F}_0,{\cal F}_{-1},{\cal F}_{-2}$ respectively.

\vskip 0,3cm

\noindent
{\bf 4.3 Infinitesimal deformations}.

\noindent
It follows from the Nijenhuis--Richardson theory that
the calculated cohomology group 
classify infinitesimal deformations of the standard embedding of
$\Vect(S^1)$ into the algebras $A(S^1)$ and ${\cal A}(S^1)$
(respectively).
One obtains the following result:

\proclaim Proposition 4.2. Every infinitesimal deformation of the
standard embedding of $\Vect(S^1)$
into $A(S^1)$ and ${\cal A}(S^1)$
is equivalent to:
\begin{equation}
f(x)\partial\mapsto f\xi+c_0f'+c_1f''\xi^{-1}+c_2f'''\xi^{-2},
\label{inf}
\end{equation}
where $c_0,c_1,c_2\in \bf R$ (or $\bf C$) are parameters.
\par

\vskip 0,3cm

To classify polynomial deformations of the standard embedding
of $\Vect(S^1)$ into ${\cal A}(S^1)$,
one needs now the integrability conditions
on parameters $c_0,c_1,c_2$.

\vskip 0,3cm

\noindent
{\bf Remark}.
We will show (cf. Section 6) that 
in the case of formal deformations
of $\Vect(S^1)$ into $A(S^1)$
the corresponding integrability conditions
are completely different.

\vskip 0,3cm

\section{Integrability condition}

\proclaim Theorem 5.1. (i) An infinitesimal deformation (\ref{inf})
corresponds to a polynomial deformation of the standard embedding
$\Vect(S^1)\hookrightarrow{\cal A}(S^1)$,
if and only if it satisfies the following condition:
\begin{equation}
6{c_0}^3c_2-3(c_0c_1)^2-18c_0c_1c_2+8{c_1}^3+9{c_2}^2=0
\label{obs}
\end{equation}
\hfill\break
(ii) The polynomial deformation corresponding to a given
infinitesimal deformation is unique up to equivalence.
\par

The nonlinear relation (\ref{obs}) is the integrability condition 
for infinitesimal deformations. Given a 1-cocycle
$C\in Z^1(\Vect(S^1);{\cal A}(S^1))$
which does not satisfy this condition, there is an obstruction
for existence of a polynomial deformation.

\vskip 0,3cm

Theorem 5.1 will be proven in the end of this section.

\vskip 0,3cm

\noindent
{\bf Remark}. 
The formula
(\ref{obs}) defines a {\it semi-cubic} parabola.
Indeed,
consider the following transformation of the parameters:
$$
\matrix{
\widetilde c_1=-2c_1+{c_0}^2\hfill\cr
\widetilde c_2=3(c_2-c_0c_1)+1\hfill\cr
}
$$
Then, the relation (\ref{obs}) is equivalent to:
$$
{\widetilde {c_1}}^3+{\widetilde {c_2}}^2=0
\eqno{(7')}
$$
\vskip 0,3cm

\noindent
{\bf 5.1 Homogeneous deformation}. 

\noindent
Consider an arbitrary polynomial deformation
of the standard embedding,
corresponding to the infinitesimal deformation (\ref{inf}):
\begin{equation}
\pi\left(f(x)\frac{d}{dx}\right)=
f\xi
+c_0f'+c_1f''\xi^{-1}+c_2f'''\xi^{-2}
+\sum_{k\in{\bf Z}}P_k(c)\pi_k(f)\xi^{-k},
\label{noh}
\end{equation}
where   $c=c_1,c_2,c_3$,
$P_k(c)$ 
are polynomials of degree $\geq2$
and $\pi_k:\Vect(S^1)\to {\cal F}_{-k}$ 
some {\it differentiable} linear maps.

Note, that since the cocycles $C_1,C_2$ and $C_3$ are defined by
differentiable maps, an arbitrary solution of the deformation problem
is also defined via differentiable maps. This follows from the
Gelfand-Fuks formalism of differentiable (or local) cohomology
(see \cite{FUK}).

\vskip 0,3cm

\noindent
{\bf Definition}.
Let us introduce a notion of homogeneity
for deformations given by differentiable maps. 
A polynomial deformation (\ref{noh}) is called {\it homogeneous} if
the sum of the degree in $\xi$
and of the order of differentiation of $f$
in each term of the right hand side
is constant.

\vskip 0,3cm

Since the cocycles $C_1,C_2$ and $C_3$ are homogeneous of order $1$,
every homogeneous deformation (\ref{noh}) corresponding to a nontrivial
infinitesimal deformation, is of homogeneity $1$:
\begin{equation}
\pi(f(x)\partial)=
f\xi
+c_0f'+c_1f''\xi^{-1}+c_2f'''\xi^{-2}
+\sum_{k\geq3}P_k(c)f^{(k+1)}\xi^{-k}
\label{hom}
\end{equation}

\proclaim Lemma 5.2. Every polynomial deformation (\ref{noh}) is equivalent to a
homogeneous deformation (\ref{hom}).
\par

\noindent
{\bf Proof}.
It is easy to see that the homomorphism equation:
$\pi([f,g])=\{\pi(f),\pi(g)\}$ preserves the homogeneity condition.
It means, that the first term in (\ref{noh}) (the term of the lowest degree
in $c$) which is not
homogeneous of degree $1$, must be a 1-cocycle.
Such a 1-cocycle is necessarily a coboundary.
Indeed, each 1-cocycle is cohomologous to a linear combination
of homogeneous of order $1$ cocycles: $C_1,C_2$ and $C_3$ (cf. Section 4.2).

The lemma follows now from the standard 
Nijenhuis-Richardson technique.
One can add (or remove) a coboundary
in any term of the polynomial deformation (\ref{noh})
to obtain an equivalent one.

Lemma 5.2 is proven.

\vskip 0,3cm

\noindent
{\bf 5.2 Uniqueness of the homogeneous deformation}.

\proclaim Proposition 5.3. Given an infinitesimal deformation (\ref{inf}),
if there exists a homogeneous polynomial deformation (\ref{hom})
corresponding to the given one, then this homogeneous polynomial deformation
is unique.
\par

\noindent
{\bf Proof}.
Substitute the formula (\ref{hom}) to the homomorphism equation.
Put
$P_0=c_0,P_1=c_1,P_2=c_2$. 
Collecting the terms
with $\xi^{-k}$
(where $k\geq3$),
one readily obtains the following identities for polynomials $P_k(c)$:

\begin{equation}
\matrix{
P_k(c)\cdot(fg'-f'g)^{(k)}=\hfill\cr\noalign{\smallskip}
\;\;\;\;\;\;\;\;\;\;\;\;\;\;\;\;\;\;\;
P_k(c)\cdot(fg^{(k+1)}+(k-1)f'g^{(k)}-
f^{(k+1)}g-(k-1)f^{(k)}g')\hfill\cr\noalign{\smallskip}
\displaystyle
\;\;\;\;\;\;\;\;\;\;\;\;\;\;\;
+\sum_{i+j=k-1}P_i(c)P_j(c)\cdot(-if^{(i+1)}g^{(j+2)}+jf^{(i+2)}g^{(j+1)}),\hfill\cr
}
\label{ide}
\end{equation}
for every $f(x), g(x)$.
Each of these identities defines a system of equations of the form:
$P_k(c)=\cdots$, where $k\geq3$ and ``$\cdots$'' means quadratic expressions
of $P_i(c)$ with $i<k$. Therefore, polynomials $P_k(c)$ with $k\geq3$
are uniquely defined by the constants $c_0,c_1$ and $c_2$.

Proposition 5.3 is proven.

\vskip 0,3cm

\noindent
{\bf 5.3 Integrability condition is necessary}.

\noindent
The first three identities (\ref{ide})
give the following system of equations:
$$
\matrix{
2P_3(c)=2c_0c_2-{c_1}^2,
\hfill\cr
5P_4(c)=3c_0P_3(c)-c_1c_2,
\hfill\cr\noalign{\smallskip}
\left\{
\matrix{
9P_5(c)=4c_0P_4(c)-c_1P_3(c)\hfill\cr
5P_5(c)=3c_1P_3(c)-2{c_2}^2\hfill\cr
}
\right.
}
$$
These equations immediately imply the relation (\ref{obs}).
This proves that this condition is necessary for integrability
of infinitesimal deformations.

\vskip 0,3cm

\noindent
{\bf 5.4 Integrability condition and the universal formula (\ref{uni})}.

\noindent
Let us prove that the condition (\ref{obs}) is sufficient
for existence of a polynomial deformation.
However, it is very difficult to solve the overdetermined
system (\ref{ide}) directly. 
We will use the formula (\ref{uni}).

The formula (\ref{uni}) obviously defines a deformation
(\ref{hom}) which is polynomial in $\lambda,\mu$.
The first two coefficients $c_0$ and $c_1$ of 
this deformation are: 
\begin{equation}
\matrix{
\displaystyle c_0=\lambda\hfill\cr\noalign{\smallskip}
\displaystyle c_1=\frac{\lambda^2}{2}-\frac{\mu^2}{2}.\hfill\cr
}
\label{par}
\end{equation}
and can be taken as {\it independent parameters}.

Fix the values of $c_0$ and $c_1$ and
consider the condition (\ref{obs}) as a quadratic equation
with $c_2$ undetermined. The two solutions can be written
(using (\ref{par})) as expressions from $\lambda$ and $\mu$:
$$
\matrix{
{c_2}^+=\frac{\lambda^3}{6}-\frac{\lambda\mu^2}{2}+
\frac{\mu^3}{3}\hfill\cr\noalign{\smallskip}
{c_2}^-=\frac{\lambda^3}{6}-\frac{\lambda\mu^2}{2}-\frac{\mu^3}{3}.\hfill\cr
}
$$
The expression ${c_2}^+$ coincides with the third coefficient
in the formula (\ref{uni}'). 
The deformation $\pi_{\lambda,\mu}$
defined by the formula (\ref{uni}), corresponds to the infinitesimal
deformation with $c_0,c_1$ given by (\ref{par})
and $c_2={c_2}^+$. 
The expression ${c_2}^-$ can be obtained from ${c_2}^+$
using the involution:
$
*:(\lambda,\mu)\mapsto(\lambda,-\mu).
$

We have shown that,
every infinitesimal deformation satisfying (\ref{obs})
corresponds to a polynomial in $c_0,c_1$ and $c_2$
deformation.
Indeed, it follows from existence of the formula (\ref{uni})
that the system (\ref{ide}) has a solution.

Theorem 5.1 is proven.

\vskip 0,3cm

\noindent
{\bf Remarks}. (a) The parameter $\mu$ is a rational parameter
on the semi-cubic parabola (\ref{obs}').
Indeed,
$\widetilde c_1=\mu^2,\widetilde c_2=\mu^3$.

(b) Suppose that 
$c_0,c_1,c_2\in\bf R$, then $\lambda$ and $\mu$ in (\ref{par})
are real if and only if ${c_0}^2\geq2c_1$.

\vskip 0,3cm

\noindent
{\bf 5.5 Proof of Theorem 1}.

We have shown that:
\hfill\break
(a) 
Every integrable infinitesimal deformation
is equivalent to (\ref{inf}) and
obeys the condition (\ref{obs}).
\hfill\break
(b)
Every polynomial deformation is equivalent to a homogeneous one.
\hfill\break
(c)
Given an infinitesimal deformation,
there exists a unique homogeneous deformation 
corresponding to the infinitesimal one.
It is given by the universal formula (\ref{uni}).

Theorem 1 is proven.





\section{Formal deformations of the embedding of $\Vect(S^1)$
into $A(S^1)$}

We classify formal deformations 
of the Lie subalgebra $\Vect(S^1)$
in $A(S^1)$.

\proclaim Theorem 2. Every formal deformation of the standard
embedding \break
$\pi:\Vect(S^1)\hookrightarrow A(S^1)$
is equivalent to one of the following deformations:
$$
\pi^t\bigg(f(x)\frac{d}{dx}\bigg)=
f\bigg(x+\frac{(1-\lambda)t}{\xi}\bigg)\xi +
\lambda tf^{\prime}\bigg(x+\frac{(1-\lambda)t}{\xi}\bigg)
\eqno{(5'')}
$$
where $\lambda\in \bf R$,
the right hand side is a (Taylor) series in $t$.
\par

In other words, there exists a one-parameter family of
formal deformations.

The explicit formula for $(\ref{uni}'')$ is:
$$
\pi^t\bigg(f(x)\frac{d}{dx}\bigg)=
\sum_{k=0}^{\infty}
(1-(k-1)\lambda)(1-\lambda)^{k-1}
\frac{t^k}{k!}
f^{(k)}(x)\xi^{-k+1}
\eqno{(5''')}
$$
note, that
$\pi^t(f(x)d\!/\! dx)=f\xi+tf'+\cdots$.

\vskip 0,3cm

Classification of
infinitesimal deformation was done in Section 4.3.
In order to prove Theorem 2, let us first classify
the infinitesimal deformations which correspond to
formal ones.

\vskip 0,3cm

\noindent
{\bf 6.1. Integrable infinitesimal deformations}.

\proclaim Proposition 6.1. An infinitesimal deformation
(\ref{inf}) correspond to a formal deformation if and only
if $c_1=c_2=0$.
\par

A 1-cocycle on $\Vect(S^1)$
corresponding to an integrable infinitesimal deformation
is, therefore, proportional to the cocycle $C_0$
from Section 4.2. 
In other words, a nontrivial formal deformation
is equivalent
to a deformation of the form:
$$
\pi^t\bigg(f(x)\frac{d}{dx}\bigg)=
f\xi+tf'+(t^2)
$$
(the constant $c_0$ in (\ref{inf}) can be chosen: $c_0=1$
up to normalization).

\vskip 0,3cm

\noindent
{\bf Proof of the proposition}. 
Consider an infinitesimal formal deformation:
$$
f(x)\frac{d}{dx}\mapsto f\xi+t(c_0f'+c_1f''\xi^{-1}+c_2f'''\xi^{-2})
$$
of the standard embedding 
$\pi:\Vect(S^1)\to A(S^1)$.
One must show that
it corresponds to a 
formal deformation $\pi^t$
if and only if $c_1=c_2=0$.

First, note that each term $\pi_k$, $k\geq1$ of a formal deformation $\pi^t$
(of the expansion (\ref{pit}))
can be chosen in a homogeneous form:
\begin{equation}
\pi_k\bigg(f(x)\frac{d}{dx}\bigg)=
\sum_{j}\alpha^k_jf^{(j+1)}\xi^j,
\label{pih}
\end{equation}
where $\alpha^k_j$ are some constants.
The proof of this fact is analogous to the one of Lemma 5.2.

Second, apply the Nijenhuis--Richardson deformation relation (\ref{deq})
(which is equivalent to the homomorphism relation
$\pi([f,g])=\{\pi(f),\pi(g)\}$).
In the same way as in Sections 5.2 and 5.3, collecting the terms with $t^k$,
one obtains the following conditions:

(a) terms with $t^2$:
$$
\matrix{
\;\;2\alpha^2_3=2c_0c_2-{c_1}^2,
\hfill\cr
\;\;5\alpha^2_4=-c_1c_2
\hfill\cr\noalign{\smallskip}
\left\{
\matrix{
\alpha^2_5=0\hfill\cr
5\alpha^2_5=-2{c_2}^2\hfill\cr
}\right.
}
$$
and therefore, $c_2=0$.

(b) terms with $t^3$:
$$
\left\{
\matrix{
9\alpha^3_5=-c_1\alpha^2_3=(1/2){c_1}^3\hfill\cr
5\alpha^3_5=3c_1\alpha^2_3=-(3/2){c_1}^3\hfill\cr
}\right.
$$
and therefore, $c_1=0$.

Proposition 6.1 is proven.

\vskip 0,3cm

\proclaim Lemma 6.2. The constants $\alpha^k_j$ in each term $\pi_k$
(given by the formula (\ref{pih}))
of the deformation $\pi^t$ satisfy the condition:
$\alpha^k_j=0$ if $j\geq k$.
\par

\noindent
{\bf Proof}.
First, one easily shows that $\alpha^k_j=0$ if $j> k$.

In the same way, using the identities (\ref{ide}), one obtains:
$\alpha_k^k=0$ for every $k\geq1$.

Lemma 6.2 is proven.

\vskip 0,3cm

For example, collecting the terms 
with $t^4$ one has
$9\alpha_5^4=4c_0\alpha_4^3-\alpha_1^2\alpha_3^2=0$ and
$\alpha_5^3=4\alpha_1^2\alpha_3^2-2(\alpha_2^2)^2=-2(\alpha_2^2)^2$,  from where
$\alpha_2^2=0$.

\vskip 0,3cm

\proclaim Lemma 6.3. Every formal deformation $\pi^t$
is equivalent to a formal deformation 
given by:
$$
\pi^t\bigg(f(x)\frac{d}{dx}\bigg)=
f\xi
+tf'+\alpha_1t^2f''\xi^{-1}+\alpha_2t^3f'''\xi^{-2}
+\sum_{k\geq3}\alpha_kt^{k+1}f^{(k+1)}\xi^{-k}
$$
where $\alpha_i$ are some constants.
\par

This means, one can take in (\ref{pih}) $\alpha^k_j=0$ 
if $j\leq k-2$.

\vskip 0,3cm

\noindent
{\bf Proof}.
Every formal deformation
is equivalent to a deformation with $\alpha^k_0=0$ in (\ref{pih}).
Indeed, constant $\alpha^k_0$ is just the coefficient behind
$t^kf'$. It can be removed (up to equivalence) 
by choosing a new formal parameter
of deformation $\widetilde t=t+t^k\alpha_0^k$.

Now, the lemma follows from Proposition 6.1 and homogeneity of the 
homomorphism condition.
Indeed, the terms with $j\leq k-2$ are independent 
and therefore, the first nonzero term 
(corresponding to the minimal value of $j$)
must be a 1-cocycle.
In the same way as in Lemma 5.2, one shows that such a 1-cocycle
is trivial and can be removed up to equivalence.

Lemma 6.3 is proven.

\vskip 0,3cm

Now, the expressions $P_k=\alpha_kt^{k+1}$ satisfy the identities
(\ref{ide}). Thus, the deformation
$\pi^t$ is given by the formula (\ref{uni}) with
$\lambda=t$. 

Theorem 2 is proven.

\section{Some properties of the main construction}

Let us study some geometric and algebraic properties of the two-parameter 
deformation (\ref{uni}).

\vskip 0,3cm

\noindent
{\bf 7.1 Deformation of $SL_2(\bf R)$-moment map}.

\noindent
Consider the standard Lie subalgebra $sl_2(\bf R)\subset\Vect(\bf R)$
generated by the vector fields:
$$
\frac{d}{dx},\;x\frac{d}{dx},\;
x^2\frac{d}{dx}.
$$
For every $\lambda$ and $\mu$, the restriction of the map $\pi_{\lambda,\mu}$ 
given by the formula (\ref{uni})
to $sl_2(\bf R)$, defines a
Hamiltonian action of $sl_2(\bf R)$ on the half-plane 
$\cal H=\{(x,\xi)\;|\xi>0\}$
endowed with the standard symplectic structure:
$\omega=dx\wedge d\xi$.
Indeed, the formal series (\ref{uni}') in this case has only
finite number of nonzero terms
and associates to each element of $sl_2(\bf R)$
a well-defined Hamiltonian function on $\cal H$.

Given a Hamiltonian action of a Lie algebra {\goth g}
on a symplectic manifold $M$, let us
recall the notion of so-called {\it moment map} from $M$
into the dual space {\goth g}$^*$ (see \cite{KIR1}).
One associates to a point $m\in M$ a linear function $\bar m$ on {\goth g}.
The definition is as follows: for every $x\in${\goth g},
$$
\langle\bar m,x\rangle:=
F_x(m),
$$
where $F_x$ is the Hamiltonian function corresponding to $x$.
If the Hamiltonian action of {\goth g} is homogeneous,
then the image of the moment map is a coadjoint orbit of {\goth g}.

In the case of $sl_2(\bf R)$, the coadjoint orbits 
on $sl_2(\bf R)^*(\simeq {\bf R}^3)$
can be identified with level surfaces of the Killing form. Explicitly,
for the coordinates on $sl_2(\bf R)^*$, dual to
the chosen generators of $sl_2(\bf R)$:
$$
y_1y_3-y_2^2={\rm const}.
$$ 
Thus, coadjoint orbits of $sl_2(\bf R)$ are cones 
(if the constant in the right hand side is zero),
one sheet of a two-sheets hyperboloid (if the constant is positive),
or a one-sheet hyperboloid (if the constant is negative).

\proclaim Proposition 7.1. The image of the half-plane $(\xi>0)$
under the $SL_2(\bf R)$-moment map is one of the following coadjoint orbits
of $sl_2(\bf R)$:
\hfill\break
(i) $\lambda=0$ or $\mu=0$, the nilpotent conic orbit;
\hfill\break
(ii) $\lambda\mu>0$, one sheet of a two-sheets hyperboloid;
\hfill\break
(i) $\lambda\mu<0$, a one-sheet hyperboloid.
\par

\noindent
{\bf Proof}. The Poisson functions corresponding to the
generators of $sl_2(\bf R)$ are:
$$
\matrix{
F_1=\xi,\hfill\cr
F_2=x\xi+\lambda,\hfill\cr
F_3=x^2\xi+2\lambda x+\lambda(\lambda-\mu)\xi^{-1},\hfill\cr
}
$$
respectively.
These functions satisfy the relation:
$F_1F_3-F_2^2=\lambda\mu$.

\vskip 0,3cm

\noindent
{\bf 7.2 The Virasoro algebra and central extension of the Lie algebra
$C^{\infty}({\bf T}^2)$}.

\noindent
Consider the Lie algebra $C^{\infty}({\bf T}^2)$
of smooth functions on the two-torus with the standard Poisson bracket.
This Lie algebra has a two-dimensional space of nontrivial central extensions:
$H^2(C^{\infty}({\bf T}^2))=H^2({\bf T}^2)={\bf R}^2$.
The corresponding 2-cocycles were defined by A.A.~Kirillov \cite{KIR} (see also
\cite{ROG}):
$$
c(F,G)=\int_{\gamma}FdG,
$$
where $F=F(x,y),G=G(x,y)$ are periodic functions:
$F(x+2\pi,y)=F(x,y+2\pi)=F(x,y)$
and $\gamma$ is a closed path.

Recall that the {\it Virasoro algebra} is the unique (up to isomorphism)
nontrivial central extension of $\Vect(S^1)$. It is given by so-called
Gelfand-Fuks cocycle:
$$
w(f(x)d\!/\! dx,g(x)d\!/\! dx)=
\int_{0}^{2\pi}f'(x)g''(x)dx
$$

Let us show how the central extensions 
of $C^{\infty}({\bf T}^2)$ are related to the Virasoro
algebra via the embedding (\ref{uni}).

\vskip 0,3cm

Let $\Vect_{\rm Pol}(S^1)$ be the Lie algebra over $\bf C$
of polynomial vector fields on $S^1$.
It is generated by:
$L_n=z^{n+1}d/dz$, where $z=e^{ix}$.
The formula (\ref{uni}) 
with $\xi=e^{iy}$ defines a family of embeddings of $\Vect_{\rm Pol}(S^1)$ into
$C^{\infty}({\bf T}^2)_{\bf C}$.

It is easy to show, that the restriction of two
basis Kirillov's cocycles
to the subalgebra 
$\Vect_{\rm Pol}(S^1)\hookrightarrow C^{\infty}({\bf T}^2)_{\bf
C}$ is proportional to the Gelfand-Fuks cocycle:
$$
\left.(\int_{\xi={\rm const}}FdG)\right|_{\Vect_{\rm Pol}(S^1)}=\lambda^2w
\;\;\;{\rm and}\;\;\;
\left.(\int_{x={\rm const}}FdG)\right|_{\Vect_{\rm Pol}(S^1)}=\lambda^2\mu^2w
$$

\vskip 1cm

\noindent
{\bf Acknowledgments}. We are grateful to F. Ziegler for fruitful
discussions. The first author would like to thank Penn. State University
for its hospitality.

\vskip 1cm


\noindent
Valentin OVSIENKO\\
{\small C.N.R.S., C.P.T.}\\
{\small  Luminy-Case 907}\\
{\small  F-13288 Marseille Cedex 9, France}\\
\\
Claude ROGER\\
{\small Institut
Girard Desargues, URA CNRS 746}
\\ {\small Universit\'e Claude Bernard - Lyon I}\\
{\small 43
bd. du 11 Novembre 1918}\\
{\small 69622 Villeurbanne Cedex, France}


\begin{thebibliography}{99}

\bibitem{AMM} {\sc F. Ammar}, {\it Syst\`emes hamiltoniens compl\`etement
integrables et d\'eformations d'alg\`ebres de Lie}.
Publications Math\'ematiques 38 (1994) 427-431.

\bibitem{FUK} {\sc D.B. Fuks}, Cohomology of infinite-dimensional Lie
algebras, Consultants Bureau, New York, 1987.

\bibitem{GER} {\sc M. Gerstenhaber}, {\it On deformations of rings and
algebras}, Annals of Math. 79, 59- .

\bibitem{KIR} {\sc A.A. Kirillov}, {\it The orbit method. Geometric
quantization}, Lectures at the University of Marilend, Preprint (1990).

\bibitem{KIR1} {\sc A.A. Kirillov}, Elements of the theory of representations,
Grundelheren der mathematische Wissenschaften 220, Springer-Verlag,
Berlin-Heidelberg-New York, 1976.

\bibitem{NR} {\sc A. Nijenhuis, R.W. Richardson}, {\it Deformations of homomorphisms
of Lie algebras}, Bull. AMS 73 (1967) 175-179.

\bibitem{OR} {\sc V. Ovsienko, C. Roger}, Deformations of Poisson brackets and
extensions of Lie algebras of contact vector fields, Russian Math. Surveys
47:6 (1992) 135-191.


\bibitem{OR1} 
{\sc  V.Yu.~Ovsienko,C.Roger}, {\it Extension of the Virasoro Group and the Virasoro Algebra by 
Modules of Tensor-Densities on $S^1$}, Funct. Anal. and its Appl.
30 (1996) No.4.

\bibitem{OR2} 
{\sc  V.Yu.~Ovsienko,C.Roger}, {\it Generalizations of Virasoro group and Virasoro algebra 
through extensions by 
modules of tensor-densities on $S^1$}, to appear in Indag. Math.

\bibitem{RIC} {\sc R.W. Richardson}, {\it , Deformations of subalgebras
of Lie algebras}, J. Diff. Geom. 3 (1969) 289-308.

\bibitem{ROG} {\sc C. Roger}, {\it Extensions centrales d'alg\`ebres et de
groupes de Lie de dimension infinie, alg\`ebre de Virasoro et g\'en\'eralisations},
Rep. on Math. Phys. 35 (1995) 225-266.

\bibitem{TAB} {\sc S. Tabachnikov},
{\it Projective connections, group Vey cocycle
and deformation quantization}, to appear in IMRN.



\end{thebibliography}
\end{document}